\documentclass[a4paper,12pt,epsfig]{article}
\usepackage[T2A]{fontenc}
\usepackage[dvips]{graphicx}
\usepackage{color}
\usepackage{graphicx}
\usepackage{amsmath,amssymb,amsthm}
\usepackage{wrapfig}

\newtheorem{lemma}{Lemma}

\newtheorem{thm}{Theorem}

\newtheorem{prop}{Proposition}[section]
\newcounter{def}

\newcommand{\dist}{\mathop{\rm dist}\nolimits}
\newcommand{\cg}{\mathop{\rm CG}\nolimits}

\begin{document}

\title{The Tammes problem for N=14}
\author{Oleg R. Musin\thanks{This research is partially supported by the NSF grant DMS-1400876 and the RFBR grant 15-01-99563.} \;  and Alexey S. Tarasov }

\date{}
\maketitle

\begin{abstract} 
The Tammes problem is to	find the arrangement of $N$ points on a unit sphere which maximizes the minimum distance between any two points. 
This  problem is presently solved  
for several values of $N$, namely for $N=3,4,6,12$ by L. Fejes T\'oth (1943); for $N=5,7,8,9$ by Sch\"utte and van der Waerden (1951); for $N=10,11$ by Danzer (1963) and for $N=24$ by Robinson (1961). Recently, we solved the Tammes problem for $N=13$. 

The optimal configuration of 14 points was conjectured more than 60 years ago. 
In the paper, we give a solution of this long-standing open problem in geometry. 
Our computer-assisted proof relies on an enumeration of the irreducible contact graphs.
\end{abstract}


\section{Introduction}

\subsection{The kissing number and thirteen spheres problem}
The {\it kissing number} $k(n)$ is the highest number of equal
non-overlapping spheres in ${\Bbb R}^n$ that touch another
sphere of the same size. 
In three dimensions the kissing number
problem asks how many white billiard balls can
{\em kiss} ({i.e.} touch) a black ball.

The most symmetrical configuration, 12 balls around
another, is achieved if the 12 balls are placed at positions
corresponding to the vertices of a regular icosahedron concentric
with the central ball. 
However, these 12 outer balls do not kiss
{one another} and {each} may be moved freely. 
This space between the balls {prompts the} question: 
{\it If you moved
all of them to one side, would a 13th ball fit?}

This problem was the subject of the famous discussion between Isaac
Newton and David Gregory in 1694. 
Most reports say that Newton
believed the answer was 12 balls, while Gregory thought that 13 might be possible. 

This problem is often called the {\it thirteen spheres problem}.
The problem was finally solved by Sch\"utte and van der Waerden
in 1953 \cite{SvdW2}. 
A subsequent two-page sketch of an elegant
proof was given by Leech \cite{Lee} in 1956.
Leech's proof was presented in the first edition of the well-known book
by Aigner and Ziegler \cite{AZ}; the authors removed this chapter from the second edition because a complete proof would have to include {too} much spherical trigonometry.

The thirteen spheres problem continues to be of interest, and
new proofs have been published in the last several years by Hsiang
\cite{Hs}, Maehara \cite{Ma, Manew} (this proof is based on Leech's proof),  B\"or\"oczky \cite{Bor},
Anstreicher \cite{Ans}, and Musin \cite{Mus13}.

Note that for $n>3$, the kissing number problem is {currently} solved only
for $n=8, 24$ \cite{Lev2,OdS}, and for $n=4$ \cite{Mus4} (see  \cite{PZ} for a beautiful exposition of this problem).

\subsection{The Tammes problem}
If $N$ unit spheres kiss the unit sphere in ${\Bbb R}^n$, then
the set of kissing points
is an arrangement on the central sphere such that the (Euclidean)
distance between any two points is at least 
1. This {observation} allows us to state the kissing
number problem in another way: {\it How many points can
be placed on the surface of ${\Bbb S}^{n-1}$ so that the angular
separation between any two points be at least $60^{\circ}$?}

{It} leads to an important generalization. 
A finite subset $X$
of ${\Bbb S}^{n-1}$ is called a {\it spherical $\psi$-code} if for
every pair  $(x,y)$ of $X$ with $x\ne y$ 
its  angular distance $\dist(x,y)$ is at least $\psi$.

\noindent
Let $X$ be a finite subset of ${\Bbb S}^{2}$. 
Denote
$$\psi(X):=\min\limits_{x,y\in X}{\{\dist(x,y)\}}, \mbox{ where } x\ne y.$$
{The set} $X$ is {then} a spherical $\psi(X)$-code.

\noindent 
Denote by $d_N$ the largest angular separation $\psi(X)$ with $|X|=N$ that can be attained
 in ${\Bbb S}^{2}$, i.e.
$$
d_N:=\max\limits_{X\subset{\Bbb S}^2}{\{\psi(X)\}}, \, \mbox{ where } \;  |X|=N.
$$
In other words,
{\it how are $N$ congruent, non-overlapping circles distributed on the sphere
when the common radius of the circles has to
be as large as possible?}

This question {is} also known as the problem of the ``inimical dictators'', namely {\it Where should $N$ dictators build their palaces on a planet so as to be as far away from each other as possible?} 
The problem
was first asked by the Dutch botanist Tammes \cite{Tam} (see \cite[Section 1.6: Problem 6]{BMP}), 
{while} examining the distribution of
openings on the pollen grains of different flowers.

The Tammes problem is presently solved 
for several values
of $N$, namely for $N=3,4,6,12$ by L. Fejes T\'oth \cite{FeT0}; for
$N=5,7,8,9$ by Sch\"utte
and van der Waerden \cite{SvdW1}; for $N=10,11$ by Danzer \cite{Dan} (for  $N=11$ see also B\"or\"oczky \cite{Bor11});
and for $N=24$ by Robinson \cite{Rob}.

\subsection{The Tammes problem for $N=13$}
The Tammes problem {for} $N=13$ is {of particular interest due to} its relation to {both} the kissing
problem and the Kepler conjecture \cite{BS,FeT}. 
Actually, this problem is equivalent to {\it the strong thirteen
spheres problem}, which asks to find the maximum radius of and an arrangement for 13 equal size non-overlapping spheres in ${\Bbb R}^3$ touching the unit sphere.

It is clear that the equality $k(3)=12$ implies $d_{13} < 60^{\circ}$. B\"or\"oczky and Szab\'o \cite{BS} proved that $d_{13} < 58.7^{\circ}$. 
Bachoc and Vallentin \cite{BV} have shown that $d_{13}
< 58.5^{\circ}$. 

We solved {the} Tammes' problem for $N=13$ in 2012 \cite{MT}. 
{We} proved that\\ {\it the arrangement $P_{13}$ of 13 points in ${\Bbb S}^2$ is the best possible, the maximal arrangement is unique up to isometry,  and   $d_{13}=\psi(P_{13})\approx 57.1367^{\circ}$.} 

In this paper, using very similar method we present a solution {to the} Tammes' problem for $N=14$.

\subsection{The Tammes problem for $N=14$}

\begin{figure}[htb]
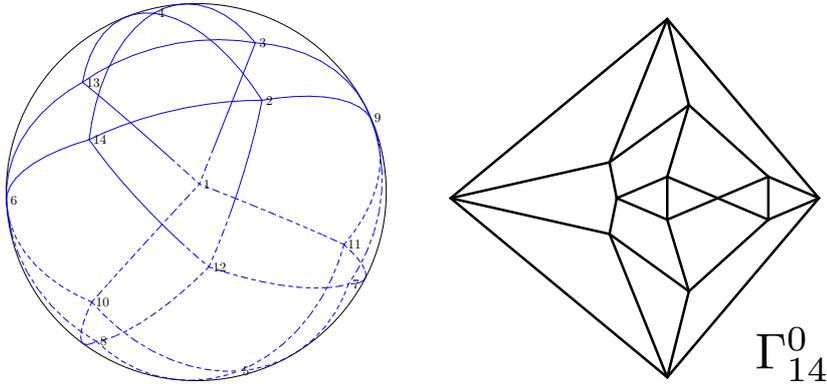

\begin{center}
\includegraphics[clip,scale=0.5]{gamma14-1.mps}
~~~~
\includegraphics[clip,scale=1.9]{tammes14-1.mps}
\end{center}
\caption{$P_{14}$ and its contact graph $\Gamma_{14}$. $\psi(P_{14})\approx 55.67057^{\circ}$.}
\label{fig1}
\end{figure}

We note that there is an arrangement of 14 points on
${\Bbb S}^{2}$ such that the distance between any two points of
the arrangement is at least  $55.67057^{\circ}$ (see \cite[Ch. VI, Sec. 4]{FeT} and http://neilsloane.com/packings/dim3/pack.3.14.txt).   
This arrangement
is shown in Fig.~\ref{fig1}.

\medskip


The first upper bound $d_{14}<58.6809^\circ$ was found in \cite{FeT0, FeT}. 
Actually, this {value} is the famous Fejes T\'oth bound
$$
d_N \le \arccos{\left(\frac{c_N}{1-c_N}\right)},  \mbox{ where } c_N:=\cos{\left(\frac{\pi N}{3N-6}\right)},
$$
for $N=14$. 

 B\"or\"oczky and Szab\'o \cite{BS14} improved the Fejes T\'oth bound and proved that $d_{14} < 58^{\circ}$. 
Bachoc and Vallentin \cite{BV} using  the SDP method have shown that $d_{14}
< 56.58^{\circ}$.

\section{Main theorem}

\begin{thm} The arrangement $P_{14}$ of $14$ points in ${\Bbb S}^2$ gives a solution of the Tammes problem,  {moreover} the maximal arrangement for $N=14$ is unique up to isometry  and   $d_{14}=\psi(P_{14})\approx 55.67057^{\circ}$.
\end{thm}

\subsection{Basic definitions}

\noindent{\bf Contact graphs.} Let $X$ be a finite set in ${\Bbb
S}^2$. 
The {\it contact graph} $\cg(X)$  is the graph with vertices
in $X$ and edges $(x,y), \, x,y\in X$ such that $\dist(x,y)=\psi(X)$.

\medskip

\noindent{\bf Shift of a single vertex.}  
We say that a vertex $x\in X$  {\it can be shifted},  if{,} in any open neighbourhood 
of $x$  there is a point  $x'\in {\Bbb S}^2$ such that $$\dist(x',X\setminus\{x\})>\dist(x,X\setminus\{x\}),$$ 
{where for a point $p\in{\Bbb S}^2$ and a finite set $Y\subset{\Bbb S}^2$  by $\dist(p,Y)$ we denote the minimum distance between $p$ and points in $Y$.}

\medskip

\begin{figure}[h]
\begin{center}
\includegraphics[clip,scale=1]{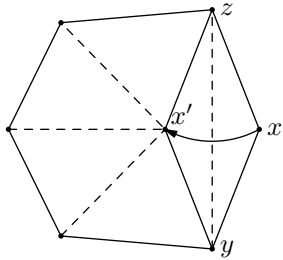}
\end{center}
\caption{Danzer's flip}
\label{fig3}
\end{figure}

\noindent{\bf Danzer's flip.}
Danzer \cite[Sec. 1]{Dan} defined the following flip. 
Let
$x,y,z$ be  vertices of $\cg(X)$ with $\dist(x,y)=\dist(x,z)=\psi(X)$. 
We say that $x$ is flipped
over $yz$ if $x$ is replaced by its mirror image $x'$ relative
to the great circle $yz$ (see Fig.~\ref{fig3}). 
We say that this flip
is {\it Danzer's flip} if $\dist(x',X\setminus\{x,y,z\}) > \psi(X)$.

\medskip

\medskip

\noindent{\bf Irreducible contact graphs.}
We say that {a} graph $\cg(X)$ is {\it irreducible} {provided it does not allow}  Danzer's flip {and no vertex in $X$ can be shifted.} 

The {concept of} irreducible contact graphs {was} invented by Sch\"utte - van der Waerden \cite{SvdW1,SvdW2}, Fejes T\'oth \cite{FeT}, and Danzer \cite{Dan}. 
Actually, in these papers as well as in our paper \cite{MT} this concept has been used for solutions of the Tammes problem with $N=7,8,9,10,11,13$.  
Recently, we enumerated all {of the} irreducible contact graphs (with and without Danzer's flip) for $N\le 11$ \cite{MT2013}.  

\medskip

\noindent{\bf $P_{14}$ and $\Gamma_{14}$.}
Denote by $P_{14}$ the arrangement of 14 points in Fig. 1. 
Let
$\Gamma_{14}:=\cg(P_{14})$. 
It is not hard  to see that the graph
$\Gamma_{14}$ is irreducible.

\medskip
\noindent{\bf Maximal graphs $G_{14}$.}
Let $X$ be a subset of ${\Bbb S}^2$ with $|X|=14$ and $\psi(X)=d_{14}$.
Denote by $G_{14}$ the graph $\cg(X)$.
Actually, this definition does not assume that $G_{14}$ is unique.
We use this designation for some $\cg(X)$ with $\psi(X)=d_{14}$.

\medskip

\begin{figure}[h]
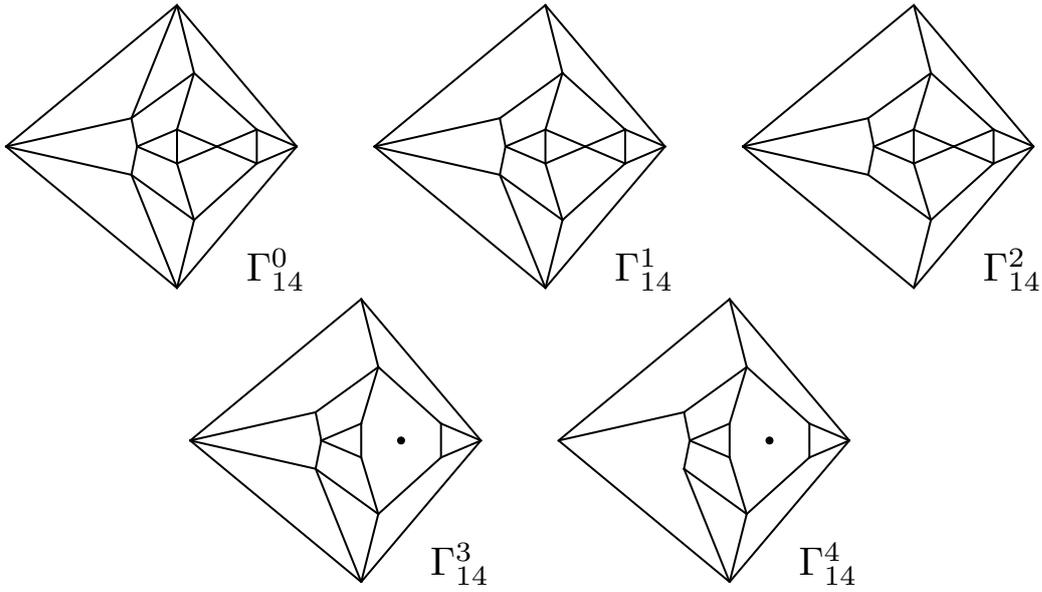

\begin{center}
\includegraphics[clip,scale=1.5]{tammes14-1.mps}
~~~~
\includegraphics[clip,scale=1.5]{tammes14-2.mps}
~~~~
\includegraphics[clip,scale=1.5]{tammes14-3.mps}
~~~~
\includegraphics[clip,scale=1.5]{tammes14-4.mps}
~~~~
\includegraphics[clip,scale=1.5]{tammes14-5.mps}
\end{center}
\caption{Graphs $\Gamma_{14}^{(i)}$ .}
\label{eliminated}
\end{figure}

\noindent{\bf Graphs $\Gamma_{14}^{(i)}$.} Let us define five
planar graphs $\Gamma_{13}^{(i)}$ (see Fig. ~\ref{eliminated}),
where $i=0,1,2,3,4$, and $\Gamma_{14}^{(0)}:=\Gamma_{14}$. 
Note
that $\Gamma_{14}^{(i)}, \, i>0$, is obtained from $\Gamma_{14}$
by removing certain edges.

\subsection{Main lemmas}

\begin{lemma} $G_{14}$ is isomorphic to  $\Gamma_{14}^{(i)}$ with $i=0,1,2,3,$ or $4$.
\end{lemma}

\begin{lemma}
 $G_{14}$ is isomorphic to $\Gamma_{14}^{(0)}$  and   $d_{14}=\psi(P_{14})$.
\end{lemma}

It is clear that Lemma 2 yields Theorem 1. Now our goal
is to prove these lemmas.

\section{Proof of Lemma 1}

Here we give a sketch of our computer {assisted} proof. For more details
see http://dcs.isa.ru/taras/tammes14/~.

\subsection{Combinatorial and geometric properties of $G_{14}$}

{The c}ombinatorial properties of $G_{N}$ {have previously been} considered {[} \cite{Dan}, \cite[Chap. VI]{FeT}, \cite{BS,BS14} and \cite{MT,MT2013} {]}.
{In particular}, for $N=14$ we have: 
\begin{prop}
{\ 
\begin{enumerate}
\item $G_{14}$ is a planar graph with 14 vertices;
\item Any vertex of $G_{14}$ is of degree $0,3,4,$ or $5$;
\item Any face of $G_{14}$ is a polygon with $3,4,5$ or $6$ vertices;
\item If  $G_{14}$ contains an isolated vertex $v$, then $v$ lies in a hexagonal face. Moreover, a hexagonal face of  $G_{14}$ cannot contain two or more isolated vertices.
\label{cor1}
\end{enumerate}
}
\end{prop}

In our papers \cite{MT,MT2013} {the} main relations 
between these parameters were considered   (\cite[Propositions 3.6--3.11]{MT}, \cite[Proposition 4.11]{MT2013}). 
Let us {list those} results {here}.  

\begin{prop}
Let $G=\cg(X)$ be an irreducible contact graph  in ${\Bbb S}^2$ with faces $F_k$. 
Let $d=\psi(X)$. 
Denote by $u_{ki}$, $i=1,\ldots,m_k$  the set of its  angles. 
Here $m_k$ denotes the number of vertices of $F_k$.
\begin{enumerate}
\item $u_{ki} <\pi$ for all $i$ and $k$.
\item $u_{ki} \geqslant \alpha(d)$ for all $i$ and $k$, where 
$$ \alpha(d):=\arccos\left( \frac{ \cos{d}}{1 + \cos{d}}\right)$$ is the angle of the equilateral spherical triangle with side length  $d$.
\item $\sum_{\tau \in I(v)} u_\tau = 2 \pi$ for all vertices $v$ of $G$. Here $I(v)$  denotes the set of    angles with the vertex in  $v$.
\item If $m_k=3$ then $F_k$ is an equilateral triangle with angles  
 $$u_{k1}=u_{k2}=u_{k3}=\alpha(d).$$
\item  In the case $m_k=4$, $F_k$ is a spherical rhombus and 
 $u_{k1}=u_{k3}$, $u_{k2}=u_{k4}$. 
 Moreover, we have the equality: 
$$
\cot{\frac{u_{k1}}{2}}\,\cot{\frac{u_{k2}}{2}} =\cos{d}.
$$

\item In the case $m_k>3$, $F_k$ is a convex equilateral spherical polygon with angles $u_{k1},\ldots,u_{m_k}$. 
{Denote  by $A_1, A_2,\ldots,A_{m_k}$ vertices of $F_k$}.
 The polygon  $F_k$ is uniquely defined (up to isometry) by  its $s:=m_k-3$ angles and  $d$. Then  functions $g_i$ and $\zeta_{ij}$, where $u_{ki}=g_i(u_{k1},\ldots,u_{ks},d)$ and $\dist(A_i,A_j)=\zeta_{i,j}(u_{k1},\ldots,u_{ks},d)$ are also uniquely defined.  
{It} follows that\\
(a) $u_{ki}=g_i(u_{k1},\ldots,u_{ks},d)$ for $i= m_k-2,m_k-1,m_k$; {and}\\
(b) $\zeta_{i,j}(u_{k1},\ldots,u_{ks},d)\geqslant d$ for all $i\ne j$.

\item  Now consider the case when  there is an isolated vertex inside  $F_k$.  
(It is only if $m_k>5$.) 
Define
$$
\lambda(u_{k1},\ldots,u_{ks},d):=\max\limits_{p\in F_k}{\min\limits_{i}\{\dist(p,A_i)\}}.
$$
Then $\lambda(u_{k1},\ldots,u_{ks},d)>d.$
\end{enumerate}
\end{prop}

\subsection{Sketch of a proof}
Our proof of Lemma 1 consists of two parts:\\
(I) Create the list $L_{14}$ of all graphs with 14 vertices that
satisfy Proposition 3.1;\\
(II) Using linear approximations and linear programming remove
from the list $L_{14}$ all graphs that do not
satisfy the {known} geometric properties of $G_{14}$ (Proposition 3.2).

\medskip

(I). To create $L_{14}$ we use the program {\it plantri}
(see \cite{PLA2}).\footnote{The authors of this program are Gunnar
Brinkmann and Brendan McKay.} 
This program  is the  isomorph-free
generator  of planar graphs, including triangulations, quadrangulations, and
convex polytopes.
({Brinkmann and McKay's} paper \cite{PLA1} describes plantri's principles of operation,
the basis for its efficiency,
and recursive algorithms behind many of its capabilities.)

The program 
generates about $1.5$ billion graphs in $L_{14}$,
i.e. graphs that satisfy {Proposition} 3.1. 
Namely, $L_{14}$
contains $409 \, 771 \,114$ graphs with triangular and quadrilateral
faces; $ 797\, 728 \, 634  $ with at least one pentagonal face  and with
triangular and quadrilaterals;
$242 \, 401 \, 773$ with at least one hexagonal face which do not contain
 isolated vertices. 

 The list of graphs with one and more isolated vertices relies on graphs in $L_N$ with $N<14$ that contain at least $14-N$ hexagons. 
For instance, the list of graphs in $L_{14}$ with exactly one isolated vertex consists of $2\,083\,967$ graphs. 
However, this list may contain isomorphic graphs.

\medskip

(II).  Let $X\subset {\Bbb S}^2$ be a finite point set such that  its contact graph $\cg(X)$ is irreducible. Properties (i)-(iv) are combinatorial properties of $\cg(X)$.  There are several geometric properties. 

Note that all faces of  $\cg(X)$ are convex. 
Since all edges of $\cg(X)$ have the same {length,} $\psi(X)$, all its faces are spherical equilateral convex polygon with number of vertices at most $\lfloor2\pi/\psi(X)\rfloor$.

 Consider now a planar graph $G$ with given faces  $\{F_k\}$ that satisfy Corollary 2.1. 
 We are going consider embeddings of this graph into  ${\Bbb S}^2$ as an irreducible contact graph  $\cg(X)$ for some $X\subset {\Bbb S}^2$.

Any embedding of $G$ in ${\Bbb S}^2$ is uniquely defined by the following list of parameters (variables):\\
(i) The edge length  $d$;\\
(ii) The set of all angles $u_{ki}$, $i=1,\ldots,m_k$ of faces $F_k$. 
Here by $m_k$ we denote the number of vertices of $F_k$.

Let us consider a graph $G$ from $L_{14}$. 
We start from
the level of approximation $\ell=1$. 
Now using Proposition 3.2  we write the 
linear equalities and inequalities {below.}

\medskip

(a)  From Proposition 3.2(3) we have 14 linear equalities  $$\sum_{\tau \in I(v)} u_\tau = 2 \pi {.}$$

(b) Since  $55.67057^\circ=0.9716\le d_{14}<0.9875=56.58^\circ$,  from Proposition 3.2(2)
we have
$$ 1.2019\le\alpha_{14}<  1.2077,$$ where $\alpha_{14}:=\alpha(d_{14})$.

\medskip

(c) For a quadrilateral $F$ with angles $u_1,u_2,u_3,u_4$ we {have} the equalities $$u_3=u_1, \; u_4=u_2,$$
and inequalities $$\alpha_{14}\le u_i\le 2\alpha_{13}, \, i=1,2.$$

\medskip

(d) For a quadrilateral $F$ we also {have} the linear inequalities 
 $$3.6057 \le u_1+u_2\le 3.7294.$$ 
These inequalities follow  from Proposition 3.2(5). 
We have $$u_2=\rho(u_1,d):=2\cot^{-1}(\tan{(u_{1}/2)}\cos d).$$ 
If we consider the maximum and minimum of  $u_1+u_2=u_1+\rho(u_1,d)$  with $u_1\in [\alpha(d),2\alpha(d)]$ and $d\ge 0.9716$, then we obtain these inequalities. 

\medskip

So from these linear equalities and inequalities we can obtain maximum and minimum values for each variable.  
It gives us a domain $D_1$ {which} contains all solutions of this system if they there exist. 
If  $D_1$ is empty, then we can remove $G$ from the list $L_{14}$.


The first step, $\ell=1$, ``kills'' almost all graphs. 
After this {first} step {all that} remained {were} $173207$ graphs without isolated vertices, $2822$ graphs
with one isolated vertex, and no graphs with two and more isolated vertices.

Next we consider $\ell=2$. In this step {$D_1$} is divided into two domains and for both we can add the same linear {constraints} as {we did when} $\ell=1$. 
Moreover, for this step we add new linear {constraints} for polygons with five and higher vertices.   

In this level we obtain the parameter
domain $D_2$. 
If this domain is empty, then $G$ cannot be embedded to ${\Bbb S}^2$ and it can be removed from $L_{14}$.  

 Actually, {for} $\ell=3$ we can repeat {the} previous step, divide  $D_2$ into two domains and obtain additional {constraints} as {we did when} $\ell=2$ for both parts independently. 

We can repeat this procedure {again and again}. 
In fact, {by} increasing  $\ell$ we increase {the} number of {sub-cases}. 
However, {in} practically 
every step some {sub-cases vanish}.

We repeat this process for $\ell=1,2,...,m$ and obtain a chain of embedded domains: $$D_m\subset\ldots\subset D_2\subset D_1.$$ 
If this chain is ended by the empty set, then $G$ can be removed from $L_{14}$.

In the case {that} a graph $G$ after 
$m$ steps still ``survive{s}'', i. e.   $D_m\ne\emptyset$,  then it is checked by numerical methods, namely by {the} so called  nonlinear ``solvers''. 
(We used, in particular, ipopt.) 
If a solution there exists, then $G$ is declared as a graph that can be embedded, and if not, then $G$ {is then removed} from $L_{14}$.

In \cite{MT2013} are given some numerical details {for} this algorithm.

\section{Proof of Lemma 2}
 In this section we present a proof of Lemma 2. 
 Actually, 
 two approaches {are considered here} -  geometric and analytic. 
 Both methods are rely on {the} geometric properties of $G_{14}$. 
 The first method we already applied to prove Lemma 2 in our solution of the Tammes problem for $N=13$ \cite{MT}. 
 In this method {by} using {the} symmetries of a graph $G=G_N$ we find certain relations between {the} variables and using {them} we prove that $\delta_N:=\psi(P_N)= d_N$. 

The geometric method is elementary, but it is not trivial and {is} relatively tricky.  
For  $G=\Gamma_{14}^{(i)}, \, i=1,2$, we found a proof that is based on the geometric approach. 
However, for the cases $G=\Gamma_{14}^{(i)}, \, i=3,4$, we could not find a simple geometric proof. 
{For those cases} we apply the analytic approach.

The idea of the analytic approach is very similar to the Connelly's ``stress matrix'' method \cite{ConRig}. 
Perhaps, this method is not {as} elementary and explicit as {the} geometric {approach}, however it works for all cases and can be {applied} with a computer assistant. 

\begin{proof}  In Section 3 we substitute all nonlinear equations by certain linear inequalities. 
Note that a statement $d_{14}\approx\delta_{14}$ is a by-product of this approximation. 
Our goal is to prove that $d_{14}=\psi(P_{14})=d(\Gamma_{14})$.

	
Lemma 1 says that $G_{14}=\Gamma_{14}^{(i)}$, where $i=0,1,2,3$, or 4. 
We are going to prove that if  $\cg(X)=\Gamma_{14}^{(i)}$ with $i>0$, then $\psi(X)<\delta_{14}$.



\medskip

\begin{figure}[htb]
\begin{center}
\includegraphics[clip,scale=1.5]{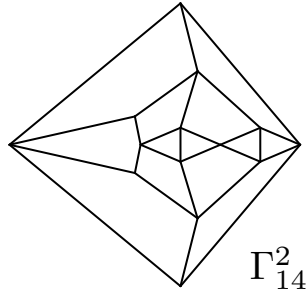}
~~~~
\end{center}
\caption{$\Gamma_{14}^{(2)}$ }
\label{gamma3}
\end{figure}

\noindent{\bf 4.1. Geometric approach: the case $\Gamma_{14}^{(i)}, \; i=1,2$.} 
Let $G=G_N$. 
Proposition 3.2(3,4,5) allows {us} to prove some equalities for variables $u_i$. 
Namely, assume that {two} vertices{,} A and B{,} of $G$ {are} adjacent {to} two triangles and two quadrilaterals. 
Then $\deg(A)=\deg(B)=4$. 
We denote the correspondent angles by  $u_i,u_j,u_0,u_0$ and $u_k,u_\ell,u_0,u_0$ respectively. 
(Here, as above,  $u_0$ denote the angle of the equilateral triangle.) 
We have 
$$
u_i+u_j+2u_0=2\pi, \, u_k+u_\ell+2u_0=2\pi, \mbox{ i. e. } u_i+u_j=u_k+u_\ell.
$$ 
If additionally, $A$ and $B$ are {the} opposite vertices of a quadrilateral $F$ in $G$, then  the equality $u_i=u_k$ (see Proposition 3.2.5)  implies the equality $u_j=u_\ell$. 
%
%


It is not hard to see that  for  the case $G=\Gamma_{14}^{(2)}$ we have  the equalities that we show in Fig.~\ref{gamma3}. 
For this graph{,} we have the following list of equations:

\begin{tabular}{ccc}
$u_0=\alpha(d)$ (i) & \\
$u_1 + u_2 + 2u_0 = 2\pi$ (ii) &
$u_{3}=\rho(u_{1},d)$ (iii) & 
$u_{4}=\rho(u_{2},d)$ (iv) \\ 
$u_5 + u_3 + 2u_0 = 2\pi$ (v) &
$u_6 + u_4 + 2u_0 = 2\pi$ (vi) &
$u_{7} =\rho(u_{5},d)$ (vii) \\
$u_{8}=\rho(u_{6},d)$ (viii) & 
$u_1 + u_{13} + 2u_7 = 2\pi$ (ix) &
$u_2 + u_{12} + 2u_8 = 2\pi$ (x) \\
$u_{15} =\rho(u_{13},d)$ (xi) &
$u_{14} =\rho(u_{12},d)$ (xii) &
$u_5 + u_{15} + u_{11} = 2\pi$ (xiii) \\
$u_6 + u_{14} + u_{11} = 2\pi$ (xiv) &
$u_0 + u_{10} + u_7+u_8 = 2\pi$ (xv) &
$u_{12} + u_{13} + 2u_{10} = 2\pi$ (xvi) 
\end{tabular}

\medskip
Let us consider two variables: $d$ and $x$, where $x:= (u_2- u_1)/2$. 
The equations (i--xvi) show that  all angles $u_i$ are uniquely {defined} by $x$ and $d$. 

From (i--viii) we have:
$$u_7 = f_7(x,d):=\rho(2\pi - 2\alpha(d)- \rho(\pi -\alpha(d) -x,d),d), $$ $$ u_8 = f_8(x,d):= \rho(2\pi - 2\alpha(d)- \rho(\pi -\alpha(d) +x,d),d).$$

Equations (ii), (ix), (x) and (xvi) yield 
$$0=  -(u_1+u_2+ 2u_0 - 2\pi) +(u_1 + u_{13}+2u_7 - 2\pi) + (u_2 +u_{12} + 2 u_8 - 2\pi)-(u_{13}+u_{12}+2u_{10} - 2\pi) $$
$$= 2 u_7 + 2 u_8 - 2u_0 -2 u_{10}$$

Therefore, $u_7 + u_8 = u_0 + u_{10}$ and (xv) implies that $u_7+u_8=\pi$. 
We obtain the equation:
$$
f_7(x,d)+f_8(x,d)=\pi. \eqno (4.1) 
$$

It is not hard to prove that if $x\in [-a,a]$, where $a>0$ is sufficiently small, then $d$ is uniquelly defined, i.e.  there is a continuous function $\theta(x)$ on $[-a,a]$ such that $d=\theta(x)$.   
Since $f_7(-x,d)=f_8(x,d)$, if $(x,d)$ is a solution of (4.1), then $(-x,d)$ is also a solution. 
It implies that the function $\theta(x)$ is even, i. e. $\theta(-x)=\theta(x)$.  
We present this function in Fig. 4.

\begin{figure}[htb]
\begin{center}
\includegraphics[clip,scale=0.5]{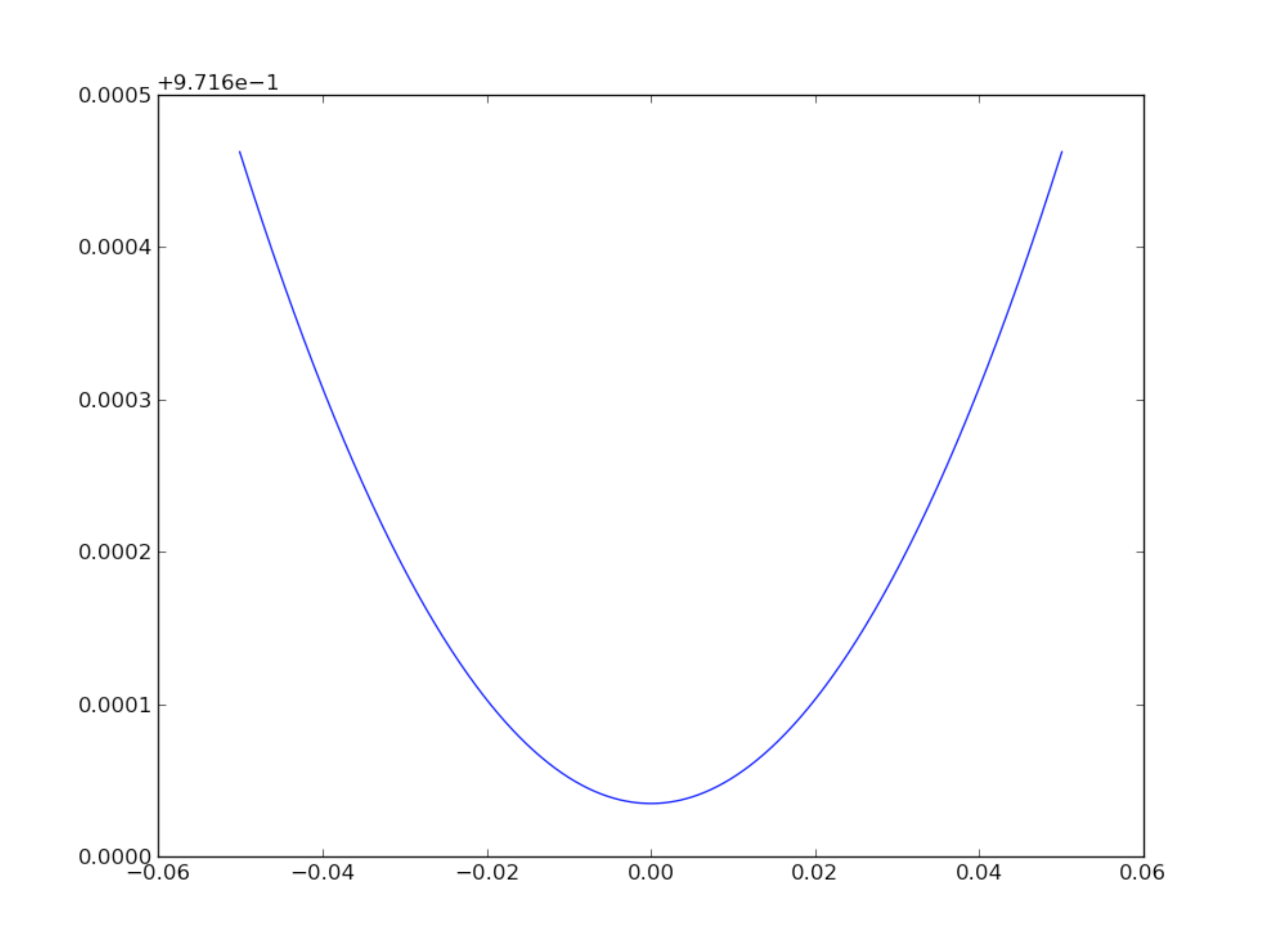}
~~~~
\end{center}
\caption{The graph of the function $d=\theta(x)$. }
\label{dbyx}
\end{figure}

Now we show that if $u_0\ge \delta_{14}$, then  $u_{13}=u_{12}=u_0$. 

From (ii) we have $u_0+u_1+x=\pi$. 
Then $u_1=\pi-x-\alpha(d)$. 
Therefore, (ix) yields 
$$
u_{13}=f_{13}(x):=2\pi - u_1-2u_7=\pi+x+\alpha(\theta(x))-2f_7(x,\theta(x)). 
$$

\begin{figure}[htb]
\begin{center}
\includegraphics[clip,scale=0.5]{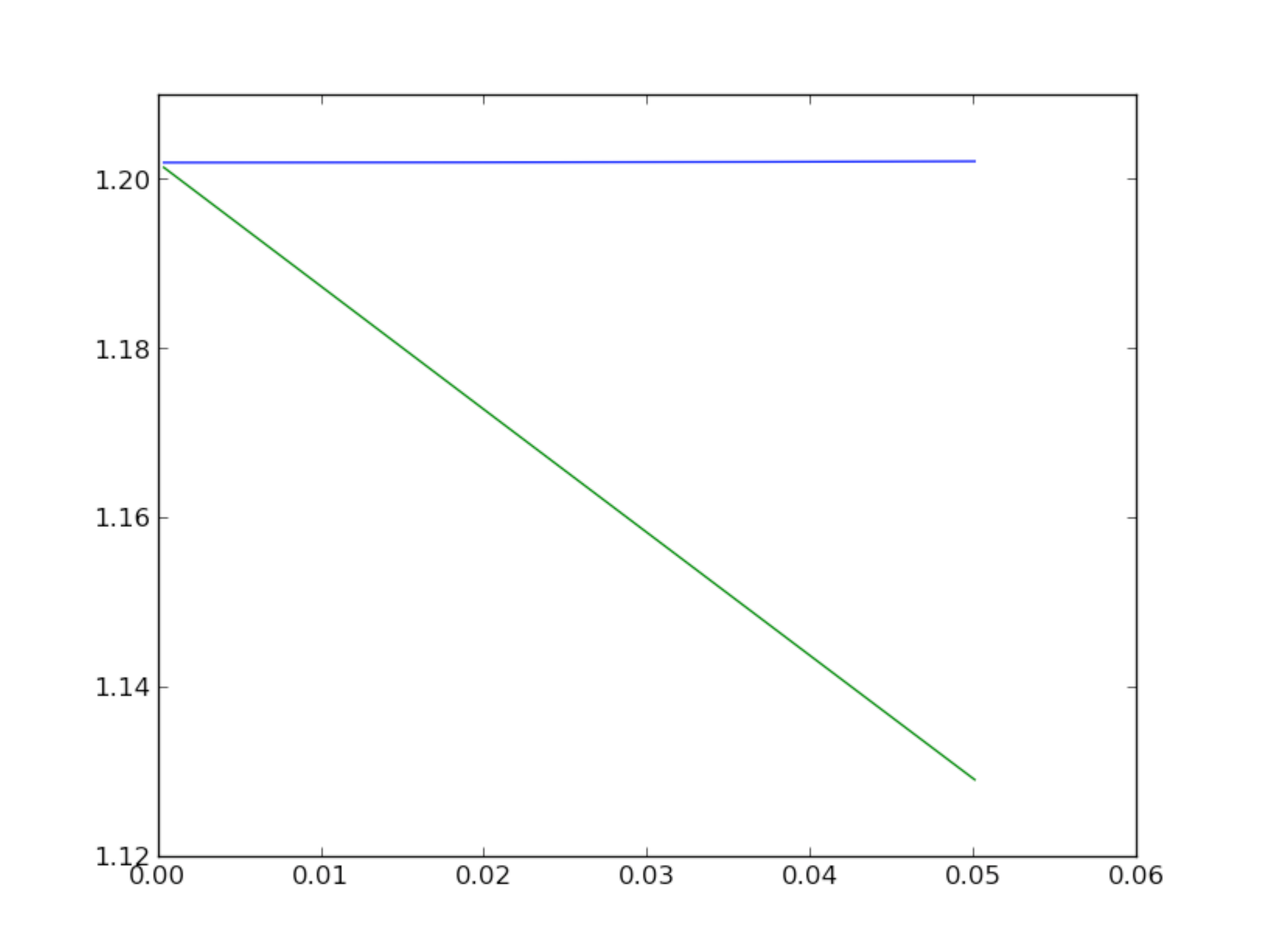}
~~~~
\end{center}
\caption{The graph of the function $u_{13}=f_{13}(x)$.}
\label{dbyx}
\end{figure}


Note that $f_{13}(0)=\delta_{14}$.  
From Fig.~\ref{dbyx} we can see that $f_{13}(x) < f_{13}(0)$ 
for $x>0$. 
Therefore,  $u_{13} < \delta_{14}\le u_0$. 
It can be rigorously proved. 
Indeed,  $f_{13}'(0)\approx -2.4587$ {which} means {that} the function $f_{13}(x)$ is {monotonically} decreasing. 
Since all $u_i\ge u_0$, we obtain $x\le 0$. 

Using the same {reasoning,} we can prove that the function $f_{12}(x)$ is {monotonically} increasing. 
Therefore, if $x<0$, then $u_{12} < u_0$. 
Thus, $u_{13}=u_{12}=u_0=\delta_{14}$. 

Note that the graph $\Gamma_{14}^{(2)}$  is a subgraph of $\Gamma_{14}^{(1)}$. 
Then $\Gamma_{14}^{(1)}$ also can not be a maximal graph {of} $G_{14}$.  

\medskip

\noindent{\bf 4.2. Analytic approach: the case $\Gamma_{14}^{(i)}, \; i=3,4$.} 

Suppose that $X=\{x_1,\ldots,x_N\}$ is a  configuration of $N$ points in the sphere ${\Bbb S}^2$ such that $\cg(X)$ is the maximal graph.  
Then points in $X$ cannot to get closer together. 
So 
any slight motion of points in $X$ cannot increase the minimal distance $\psi(X)$. 
Therefore, $X$
is an infinitesimally rigid configuration \cite{ConRig}.

We say that an $N\times N$ symmetric matrix $(\omega_{ij})$ is the {\it equilibrium stress} matrix if to each pair of distinct vertices $\{i,j\}$ of $G_N:=\cg(X)$ we have $\omega_{ij} \ge 0$,  $\omega_{ij}=0$ when  $\{i,j\}$ is not an edge of $G_N$, and for each $i$, the equilibrium equation 
$$ F_i: = \sum\limits_{j,\, j\ne i} {\omega_{ij} e_{ij}} = 0 \eqno (4.2)$$
holds. 
Here $e_{ij}$ is the unit tangent  vector at the point $x_i$ to the great circle {that} passes {through} the points $x_i$ and $x_j$.

Actually, these conditions for the equilibrium stress matrix can be derived from the  Karush-Kuhn-Tucker conditions \cite[p. 244]{BV}, where stresses $\omega_{ij}$  correspond to Lagrange multipliers and the equilibrium equation 
{corresponds} to the stationary condition. 
Note that,  the inequality $\omega_{ij} \ge 0$ holds because  $\dist(x_i,x_j)=\psi(X)$ cannot be increased.

Now using conditions (4.2) we set up a system of linear inequalities. 
Since the linear programming shows that this system  has no solution, we obtain that $\Gamma_{14}^{(i)}, \; i=3,4$, cannot be maximal graphs. 

Denote by $J(i)$ the set of all indexes $j:\omega_{ij}\ne0$, i. e. $\omega_{ij}>0$. 
We have 
$$
\sum\limits_{j\in J(i)}{\omega_{ij}e_{ij}}=0. \eqno (4.3)
$$

In fact, after computations we have approximation intervals for all parameters  of $\Gamma_{14}^{(i)}.$ Therefore, we can compute intervals also for all $c_{ij}$ and $s_{ij}$, where  
$$e_{ij}=(c_{ij},s_{ij})\; \mbox{ and } \; c_{ij}^2+s_{ij}^2=1.$$ 
Let $c_{ij}\in[p_{ij},q_{ij}]$ and $s_{ij}\in[u_{ij},v_{ij}]$. 
Here $q_{ij}-p_{ij}$ and $v_{ij}-u_{ij}$ are sufficiently small numbers. 
Thus, (4.3) implies  
$$
\sum\limits_{j\in J(i)}{\omega_{ij}p_{ij}}\le0, \; \; \sum\limits_{j\in J(i)}{\omega_{ij}u_{ij}}\le0, \; \;
\sum\limits_{j\in J(i)}{\omega_{ij}q_{ij}}\ge0, \; \; \text{{and} } \sum\limits_{j\in J(i)}{\omega_{ij}v_{ij}}\ge0. \eqno (4.4) 
$$

{Now let} us add to (4.4) the normalization condition: 
$$\sum\limits_{i,j} \omega_{ij} = 1.$$
As we mentioned above this system has no solution. 
Thus, 
 $\Gamma_{14}^{(3)}$ and $\Gamma_{14}^{(4)}$ cannot be maximal graphs. 

\medskip



\end{proof}

\medskip

{\bf Acknowledgment.} We  wish to thank James Maissen for his helpful comments and corrections 
on this written work.

\medskip

\medskip

\medskip

\medskip

\medskip

\medskip

O. R. Musin, University of Texas at
Brownsville.

 {\it E-mail:} oleg.musin@utb.edu

\medskip

A. S. Tarasov, IITP RAS

{\it E-mail:} tarasov.alexey@gmail.com

\end{document}